\documentclass[12pt, reqno]{amsart}
\usepackage{epsfig}
\usepackage{enumitem}
\usepackage{multicol}
\usepackage{makecell}
\def\footnoterule{\kern-3pt \hrule width 1.9 cm \kern2.6pt}
\begin{document}
	% Type the title of the article here
	\centerline{\bf \small{Polygroupoid, Polyquasigroup, Polyloop and their  Nuclei}} 
	\vskip 0.4 cm 
	%Type the name(s) of the author(s)
	\centerline{\scriptsize{K. G. ILORI$^{2}$, T. G. JAIY\'E\d{O}L\'A\footnote{Corresponding 
				author}, O.O. OYEBOLA$^{3}$, O. B. OGUNFOLU$^{4}$ AND E. A. ALHASSAN$^{5}$}}
	\vskip 0.4 cm
	\begin{center}
		\begin{minipage}[b]{9.0cm}
			\scriptsize{
				\noindent % Insert the abstract of your paper here
				ABSTRACT. In this paper, new hyper-algebraic structures called polygroupoid, polyquasigroup and polyloop were introduced with concrete examples given. The first, second, third and fourth left (middle, right) nuclei of polygroupoid were introduced and studied. It was shown that first left (middle, right) nuclei of a polygroupoid is contained in second, third and fourth left (middle, right) nuclei of the polygroupoid. Hence, the second, third and fourth left (middle, right) nuclei of a polygroupoid generalize the first left (middle, right) nuclei  of the polygroupoid. Examples of polyquasigroups and polyloops that satisfy the above results were provided. 
			}
		\end{minipage}
	\end{center}
	\vskip 0.2 cm
	\noindent % Insert a maximum of 5 Keywords and phrases 
	\par
	\noindent % Insert the AMS classification
	\par
	\vskip 0.5 cm
	{\small{
			\begin{center}{\section{INTRODUCTION}}\end{center}}
	}
	\vskip 0.3 cm
	\noindent %Insert the content of your introduction
\noindent
{\bf 1.1}	Hypergroup, $H_{v}$-group and Polygroup
	\vskip 0.3 cm
The hypergroup notion was introduced in 1934 by F. Marty \cite{ilo1}, at the
$8^{th}$ Congress of Scandinavian Mathematicians. Hypergroups are a suitable generalization
of groups. We know in a group, the composition of two elements is an
element, while in a hypergroup, the composition of two elements is a set. The motivation example was the following: Let $G$ be a group and $H$ be any subgroup of $G$. Then, $G/H=\{xH~|~x \in G\}$ becomes a hypergroup where the hyperoperation is defined in a usual manner:
\begin{displaymath}
	aH\circ bH=\{cH~|~c \in aH\cdot bH\}~\textrm{for all}~a,b\in G.
\end{displaymath}

Several kinds of hypergroups have been intensively studied, such as regular hypergroups, reversible regular hypergroups, canonical hypergroups,
cogroups and cyclic hypergroups. The situations that occur in hypergroup
theory, are often extremely diversiﬁed and complex with respect to group
theory.
	\vskip 0.3 cm
\noindent
{\bf Definition 1.1:} (Semihypergroup, Quasihypergroup, Hypergroup, $H_v$-group)
\noindent	
	\vskip 0.3 cm
\noindent 
An hypergroupoid $(H,\circ)$ is the pair of a non-empty set $H$ with an hyperoperation $\circ :H\times H\to P(H)\backslash\{\emptyset\}$ defined on it. An hypergroupoid $(H,\circ)$ is called a semihypergroup if
	\begin{description}
		\item[(i)] it obeys the associativity law $a \circ (b \circ c) = (a \circ b) \circ c$ for all $a, b, c\in H$, which means that
		\begin{displaymath}
			\displaystyle \bigcup_{u\in~a \circ b}u \circ c = \bigcup_{v\in~b \circ c}a \circ v
		\end{displaymath}
\noindent 
An hypergroupoid $(H,\circ)$ is called a quasihypergroup if
		\item[(ii)] it obeys the reproduction axiom $x\circ H=H=H\circ x$ for all $x\in H$.
	\end{description}
	An hypergroupoid $(H,\circ)$ is called an hypergroup if it is a semihypergroup and a quasihypergroup.
	
	A hypergroupoid $(H,\circ)$ is called an $H_v$-group if it is a quasihypergroup and it obeys the weak associativity (WASS) condition
	\begin{description}
		\item[(iii)] $x \circ (y \circ z) \cap (x \circ y) \circ z \neq \emptyset$ for all $a, b, c\in H$.
	\end{description}
%\end{mydef}
The second condition is frequently used in the form: Given $a$, $b\in H$, there exist $x, y\in~H$ such that $b\in~a \circ x$ and $b\in~y \circ a$.
Hence, an hypergroup (of Marty) is equivalent to a multigroup of Dresher and Ore \cite{dresh}.

\noindent
{\bf Definition 1.2:}(Polygroup, Davvaz \cite{tolu1})
	
	A polygroup is a system $\wp = < P, \cdot, e, {}^{-1} >$, where $e \in P$, $^{-1}$ is a unitary operation on $P$, `$\cdot$' maps $P\times P$ into the non-empty subsets of $P$, and the following axioms hold for all $x, y, z \in P$:
			\begin{description}
			\item[(P1)]$(x \cdot y) \cdot z = x \cdot (y \cdot z)$,
			\item[(P2)]$e \cdot x = x \cdot e = x$
			\item[(P3)]$x \in y \cdot z$ implies $y \in x \cdot z^{-1}$ and $z \in y^{-1} \cdot x$.
		\end{description}

\noindent
{\bf Remark 1.1:}
	A polygroup is a special type of hypergroup.

\noindent
{\bf Example 1.1:} (Double Coset Algebra, Davvaz \cite{tolu1})
	
	Suppose that $H$ is a subgroup of a group $G$. Define a system $G//H = \big< \{HgH~|~g \in G\}, \ast, H, ^{-I} \big>$, where $(H g H)^{-I} = H g^{-1} H$ and
	\begin{displaymath}
		(Hg_{1}H) \ast (Hg_{2}H) = \{Hg_{1}hg_{2}H~|~h \in H\}.
	\end{displaymath}
	The algebra of double cosets $G//H$ is a polygroup.

Articles and books have been published up-till date on hyperstructures, among which are Corsini and Leoreanu \cite{corsini}, Vougiouklis \cite{voug} and Davvaz \cite{tolu1,New2}, Davvaz and Vougiouklis \cite{New1}. Many of these were dedicated to the theories and applications of hyperstructures in other topics.  As it was stated in \cite{tolu1}, some of the fields connected with hyperstructures are Geometry, Codes, Cryptography and Probability, Automata, Artiﬁcial Intelligence, Median Algebras, Relation Algebras, C$^*$-algebras, Boolean Algebras, Categories, Topology, Binary Relations, Graphs and Hypergraphs, Lattices and Hyperlattices, Fuzzy Sets and Rough Sets, Intuitionistic  Fuzzy  Hyperalgebras, Generalized Dynamical Systems and so on.

Oyebola and Jaiy\'e\d{o}l\'a \cite{Oye1} investigated non-associative properties in algebraic hyperstructures as it plays out in the biological inheritance which is expressed in the genotypic and
phenotypic information that are passed to the progenies from the parental traits. These structures were found to be hypergroupoids or hyperquasigroups which obey $1$-variable
identity ($3$-power associativity) or $2$-variable identities (LAP, RAP or flexibility) or $3$-variable
identities (extra-$1$ or extra-$2$ or extra-$3$). Such hyperstructures were termed to be $3$-power
associative, flexible, left (right) alternative or extra; in their precise measure of weakness in
associativity. This was with the objective to valuate with precision the non-associativity of weak associative properties (WASS) in algebraic structures derived from some biological inheritance crossing. Ilori et al. \cite{chem1} carried out the analysis of weak hyper-algebraic structures that represent dismutation reaction.

In this current paper, new hyper-algebraic structures called polygroupoid, polyquasigroup and polyloop will be introduced and concrete examples given. The first, second, third and fourth left (middle, right) nuclei of polygroupoid will be introduced and studied. Examples of polygroupoid, polyquasigroups and polyloops that satisfy the above results will be provided. 
	\vskip 0.3 cm
	
	\noindent
	\begin{center}
		\section{PRELIMINARY}
	\end{center}
	\subsection{Quasigroup and Loop}
	\paragraph{}
	\vskip 0.3 cm
	Let $G$  be a non-empty set. Define a binary operation ($\cdot$) on
	$G$. $(G,\cdot)$ is called a groupoid if $G$ is closed under the
	binary operation ($\cdot$). A groupoid $(G,\cdot)$ is called a
	quasigroup if the equations $a\cdot x=b$ and $y\cdot c=d$ have
	unique solutions for $x$ and $y$ for all $a,b,c,d\in G$. A
	quasigroup $(G,\cdot)$ is called a loop if there exists a unique
	element $e\in G$ called the identity element such that $x\cdot
	e = e\cdot x = x$ for all $x\in G$.
	
	Alimpic \cite{alim} considered a generalization of the nucleus of a groupoid on the $n$-ary case, it was established that $(N_{\lambda},\cdot )$, $(N_{\mu}, \cdot)$ and $(N_{\rho}, \cdot)$ form semigroup.
	
\noindent
{\bf Definition 2.1:}
		A \textit{quasigroup} $(G,\cdot ,\backslash ,/)$ is a set $G$ together
		with three binary operations ($\cdot $), ($/$), ($\backslash$) such that
		\begin{description}
			\item[(i)] $x\cdot (x\backslash y)=y$, $(y/x)\cdot x=y$ for all
			$x,y\in G$,
			\item[(ii)] $x\backslash (x\cdot y)=y$, $(y\cdot x)/x=y$ for all
			$x,y\in G$.
		\end{description}
		Let $(G,\cdot ,\backslash ,/)$ be a \textit{quasigroup}. Then, $(G,\cdot ,\backslash ,/,e)$ is called a \textit{loop} if there is a nullary operation $e$ such that
		\begin{description} 
			\item[(iii)] $x\backslash x=y/y$ or $e\cdot x=x$ for all
			$x,y\in G$.
		\end{description}

	We also stipulate that ($/$) and ($\backslash$) have higher priority
	than ($\cdot $) among factors to be multiplied. For instance,
	$x\cdot y/z$ and $x\cdot y\backslash z$ stand for $x(y/z)$ and
	$x(y\backslash z)$ respectively.
	
	\noindent
	{\bf Definition 2.2:}\label{nucl1} (Nuclei of groupoid)
		
		Let $(G,\cdot)$ be a groupoid. The left nucleus $N_{\lambda}$ (middle nucleus $N_{\mu}$, right nucleus $N_{\rho}$) is the set of all left (middle, right) nuclear elements in $(G,\cdot)$ and the nucleus $N$ of $(G, \cdot)$ is given by $N=N_{\lambda}\cap N_{\mu}\cap N_{\rho}$.
		\begin{center}
			$N_{\lambda}(G, \cdot) = \{a \in G ~\mid~ a \cdot (x \cdot y) = (a \cdot x) \cdot y ~\forall~ x, y \in G\}$
			
			$N_{\mu}(G, \cdot) = \{a \in G ~\mid~ (x \cdot a) \cdot y = x \cdot (a \cdot y) ~\forall~ x, y \in G\}$
			
			$N_{\rho}(G, \cdot) = \{a \in G ~\mid~ (x \cdot y) \cdot a = x \cdot (y \cdot a) ~\forall~ x, y \in G\}$
		\end{center}

	Readers can consult the books \cite{Bruck b,  davidref:10, Pflugfelder,springerbook} for a general background in loop theory.
	
	\subsection{Polyquasigroup and Polyloop}
	\paragraph{}
	We now introduce new hyper-algebraic structures (non-associative) and some notions with which they will be studied.
	\paragraph{}
	\noindent
	{\bf Definition 2.3:} (Hypergroupoid, Hyperquasigroup and Hyperloop)
		
		Let $H$ be a non-empty set and $\cdot : H \times H \rightarrow~ \mathcal{P}^{\ast}(H)$ be an hyperoperation, where $\mathcal{P}^{\ast}(H)$ is the family of non-empty subsets of $H$. The couple $(H, \cdot)$ is called a hypergroupoid (or polygroupoid). A hypergroupoid will be called an hyperquasigroup, if it satisfies the reproduction axiom. i.e. for any $a \in H$, then
		\begin{center}
			$a \cdot H = H \cdot a = H$
		\end{center}If in addition, there exist $e\in H$ such that $x\in x \cdot e = e \cdot x$ for all $x\in H$, then $(H, \circ, e)$ is called an hyperloop. Thus, an associative hyperquasigroup is an hypergroup.

	\noindent
	{\bf Remark 2.1:}
		Note that Usan and Galic \cite{usan} used the term `hyperquasigroup' in a different way; an hyperquasigroup is an hypergroupoid with division (reproductive axiom and hypercancellation).

	\noindent
	{\bf Definition 2.4:}(Geometry Hyperquasigroup, Tallini \cite{geometric})
		\begin{enumerate}
			\item
			Let $(H,\cdot )$ be an hypergroupoid. If the following axioms hold: 
			\begin{enumerate}
				\item for all $h\in H,~h\cdot h = h$ (Tallini 1 axiom);
				\item for all $h_1,h_2\in H,~h_1,h_2\in h_1 \cdot h_2$ (Tallini 2 axiom); 
				\item for all $h_1,h_2\in H$ and  for all $g_1,g_2\in H,~g_1\ne g_2\implies h_1\cdot h_2\subseteq g_1\cdot g_2$ (Tallini 3 axiom);
			\end{enumerate}
			then, $(H,\cdot )$ is called a geometric hyperquasigroup. 
			\item In an hypergroupoid $(H,\cdot )$, the axiom
			\begin{displaymath}\forall~h_1,h_2\in H,~\forall~h_3\in h_1\cdot h_2\implies h_1\cdot (h_2\cdot h_3)=(h_1\cdot h_2)\cdot h_3=h_1\cdot h_2
			\end{displaymath}
			will be called the Tallini 4 axiom.
			\item
			Let ${\mathcal{H}}=\{x\cdot y : (x, y)\in H^2\}\subset H^2$. $(H,{\mathcal{H}})$ is called a line space associated with the hypergroupoid $(H,\cdot )$ if there is a unique element in ${\mathcal{H}}$ through any two distinct elements in $H$. 
		\end{enumerate}

	\noindent
	{\bf Remark 2.2:}
		Note that in an hypergroupoid $(H,\cdot )$, Tallini 2 and Tallini 3 axioms implies the property: 
		\begin{gather*}
			\textrm{for all $h_1,h_2\in H$ and  for all $g_1,g_2\in H,~g_1\ne g_2\implies$}\\
				\textrm{$h_1\cdot h_2= g_1\cdot g_2$--Tallini 5 axiom}
		\end{gather*}
		In fact, axioms Tallini 1, Tallini 2 and Tallini 5 imply commutativity and Tallini 4 axiom.

	\noindent
	{\bf Lemma 2.1:}(Tallini \cite{geometric})	
		\begin{enumerate}
			\item The subhyperquasigroups of a hyperquasigroup  $(H,\cdot )$ form a closure system  coincident with the closure system of the subspaces in the line space associated with $(H,\cdot )$
			\item A geometric hyperquasigroup $(H,\cdot )$ is a hypergroup if and only if $(H,{\mathcal{H}})$  is a projective space. 
		\end{enumerate}

	\noindent
	{\bf Definition 2.5:}\label{kendef2} (Polygroupoid, Polyquasigroup, Polyloop, Multiloop)
		
		Let $\mathcal{M} = (P, \cdot)$ be a polygroupoid. Let $e \in P$ and $/ : P \times P \rightarrow \mathcal{P}^{\ast}(H)$ and $\backslash : P \times P \rightarrow \mathcal{P}^{\ast}(H)$ such that
		\begin{description}
			\item[(a)] (i) $x \in (x \cdot y)/ y$ (ii) $x\in (x/ y) \cdot y$ (iii) $x \in y\backslash (y \cdot x)$ (iv) $x\in y \cdot (y\backslash x)$ for all $x, y \in P$, then $(P, \cdot, \backslash,/)$ will be called a polyquasigroup.
			\item[(b)] $x \cdot e = e \cdot x = x$ for all $x \in P$ and $(P, \cdot,\backslash, /)$ is a polyquasigroup. Then $(P, \cdot,\backslash, /, e)$ will be called a polyloop.
			\item[(c)] $x\in x \cdot e = e \cdot x$ for all $x \in P$ and $(P, \cdot,\backslash, /)$ is a polyquasigroup. Then $(P, \cdot,\backslash, /, e)$ will be called a multiloop.
			\item[(d)] $(x \cdot y) \cdot z = x \cdot (y \cdot z)$ for all $x, y, z \in P$ and $(P, \cdot,\backslash, /)$ is a polyloop. Then $(P, \cdot,\backslash, /)$ will be called an associative polyloop.
		\end{description}

	\noindent
	{\bf Remark 2.3:}\label{kenrem1} In a polygroupoid $(P, \cdot)$, we shall often write $x \cdot y$ as $xy~\forall~x, y \in P$. In Definition~\ref{kendef2}, (a)(i) means $b\in x \cdot a \Rightarrow x\in P$; (a)(ii) means $b\in x/ a~\Rightarrow~x \in P$; (a)(iii) means $b\in a \cdot y~\Rightarrow~y \in P$; (a)(iv) means $b\in a\backslash y~\Rightarrow~y\in P$. Thus, $b\in x\backslash a ~\Rightarrow~x\in P$ and $b\in a / y~\Rightarrow~y\in P$ for any two $a, b\in P$.

	\par
	\noindent
	\vskip 0.5 cm
	{\small{
			\begin{center}{\section{MAIN RESULTS}}\end{center} }}
	\vskip 0.3 cm
	\noindent
	\par			%paragraph in the article
	We first give some examples of  polygroupoid, polyquasigroups and polyloops that are neither hypergroups nor polygroups. They are also not geometric hyperquasigroups.
	\par
	\noindent				
	{\bf Example 3.1:} (Example of polygroupoid)
	
	Let $G = \{1, 2, 3, 4, 5, 6\}$ and define the hyperoperation $\odot: G \times G \rightarrow \mathcal{P}^{*}(G)$ as in Table~\ref{newexampolygroupoid3}. Then, $(G, \odot)$ is a polygroupoid. It is however neither a polyquasigroup nor a polyloop.
	\begin{table}[!hbp]
		\begin{center}\caption{Multiplication table for polygroupoid $(G, \odot)$}\label{newexampolygroupoid3}
			\begin{tabular}{|c|c|c|c|c|c|c|}
				\hline
				$\odot$	& 1 & 2 & 3 & 4 & 5 & 6 \\
				\hline
				1	& \{1\} & \{3, 4\} & \{1, 5\} & \{1, 4\} & \{1, 3\} & \{4, 5\}   \\
				\hline
				2	& \{3, 4\} & \{6\} & \{2, 3\} & \{4, 6\} & \{3, 6\} & \{2, 6\}  \\
				\hline
				3	& \{1, 5\} & \{2, 3\} & \{5\} & \{1, 2\} & \{3, 5\} & \{2, 5\} \\
				\hline
				4	& \{1, 4\} & \{4, 6\} & \{1, 2\} & \{4\} & \{1, 6\} & \{2, 4\}  \\
				\hline
				5	& \{1, 3\} & \{3, 6\} & \{3, 5\} & \{1, 6\} & \{3\} & \{5, 6\}  \\
				\hline
				6   & \{4, 5\} & \{2, 6\} & \{2, 5\} & \{2, 4\} & \{5, 6\} & \{2\} \\
				\hline
			\end{tabular}
		\end{center}
	\end{table}
	
		\par
	\noindent				
	{\bf Example 3.2:} (Example of polyquasigroup)
	
	Let $G = \{1, 2, 3, 4, 5, 6, 7\}$ and define the hyperoperations $\odot,\nwarrow,\nearrow: G \times G \rightarrow \mathcal{P}^{*}(G)$ as in 
	Table \ref{newexampolyloop3}, Table \ref{quasiexampolyloop4}, Table \ref{quasiexampolyloop3}. Then, $(G, \odot, \nwarrow, \nearrow)$, $(G, \nwarrow, \odot, \nwarrow)$ and  $(G, \nearrow, \nearrow,\odot)$ are polyquasigroups. They are however not polyloops.
	\begin{table}[!hbp]
		\begin{center}\caption{Multiplication table for  polyquasigroup $(G, \odot, \nwarrow, \nearrow)$}\label{newexampolyloop3}
			\begin{tabular}{|c|c|c|c|c|c|c|c|}
				\hline
				$\odot$	& 1 & 2 & 3 & 4 & 5 & 6 & 7 \\
				\hline
				1	& 2 & \{4,5\} & \{1, 6\} & \{1, 4\} & \{3, 7\} & \{3,6\} & \{6, 7\} \\
				\hline
				2	& \{4, 5\} & \{7\} & \{2, 3\} & \{4, 6\} & \{2, 6\} & \{1, 5\} & \{1, 3\} \\
				\hline
				3	& \{1, 6\} & \{2, 3\} & \{5\} & \{1, 7\} & \{2, 4\} & \{4, 7\} & \{3, 6\} \\
				\hline
				4	& \{1, 4\} & \{4, 6\} & \{1, 7\} & \{3\} & \{5, 6\} & \{2, 7\} & \{2, 5\} \\
				\hline
				5	& \{3, 7\} & \{2, 6\} & \{2, 4\} & \{5, 6\} & \{1\} & \{3, 4\} & \{5, 7\} \\
				\hline
				6	& \{3, 6\} & \{1, 5\} & \{4, 7\} & \{2, 7\} & \{3, 4\} & \{6\} & \{1, 2\} \\
				\hline
				7	& \{6, 7\} & \{1, 3\} & \{3, 6\} & \{2, 5\} & \{5, 7\} & \{1, 2\} & \{4\} \\
				\hline
			\end{tabular}
		\end{center}
	\end{table}
	
	\begin{table}[!h]
		\begin{center}\caption{Multiplication table for  polyquasigroup $(G, \nwarrow, \odot, \nwarrow)$}\label{quasiexampolyloop4}
			\begin{tabular}{|c|c|c|c|c|c|c|c|}
				\hline
				$\nwarrow$	& 1 & 2 & 3 & 4 & 5 & 6 & 7 \\
				\hline
				1	& \{3, 4\} & \{6, 7\} & \{1, 4\} & \{1, 3\} & \{5\} & \{2,7\} & \{2, 6\} \\
				\hline
				2	& \{1\} & \{3, 5\} & \{2, 5\} & \{6, 7\} & \{2, 3\} & \{4, 7\} & \{4, 6\} \\
				\hline
				3	& \{5, 6\} & \{3, 7\} & \{2, 7\} & \{4\} & \{1, 6\} & \{1, 5\} & \{2, 3\} \\
				\hline
				4	& \{2, 4\} & \{1, 4\} & \{5, 6\} & \{1, 2\} & \{3, 6\} & \{3, 5\} & \{7\} \\
				\hline
				5	& \{2, 6\} & \{1, 6\} & \{3\} & \{5, 7\} & \{4, 7\} & \{1, 2\} & \{4, 5\} \\
				\hline
				6	& \{3, 7\} & \{4, 5\} & \{1, 7\} & \{2, 5\} & \{2, 4\} & \{6\} & \{1, 3\} \\
				\hline
				7	& \{5, 7\} & \{2\} & \{4, 6\} & \{3, 6\} & \{1, 7\} & \{3, 4\} & \{1, 5\} \\
				\hline
			\end{tabular}
		\end{center}
	\end{table}
	
	\begin{table}[!h]
		\begin{center}\caption{Multiplication table for polyquasigroup $(G, \nearrow, \nearrow ,\odot)$}\label{quasiexampolyloop3}
			\begin{tabular}{|c|c|c|c|c|c|c|c|}
				\hline
				$\nearrow$	& 1 & 2 & 3 & 4 & 5 & 6 & 7 \\
				\hline
				1	& \{3, 4\} & \{1\} & \{5, 6\} & \{2, 4\} & \{2, 6\} & \{3,7\} & \{5, 7\} \\
				\hline
				2	& \{6, 7\} & \{3, 5\} & \{3, 7\} & \{1, 4\} & \{1, 6\} & \{4, 5\} & \{2\} \\
				\hline
				3	& \{1, 4\} & \{2, 5\} & \{2, 7\} & \{5, 6\} & \{3\} & \{1, 7\} & \{4, 6\} \\
				\hline
				4	& \{1, 3\} & \{6, 7\} & \{4\} & \{1, 2\} & \{5, 7\} & \{2, 5\} & \{3, 6\} \\
				\hline
				5	& \{5\} & \{2, 3\} & \{1, 6\} & \{3, 6\} & \{4, 7\} & \{2, 4\} & \{1, 7\} \\
				\hline
				6	& \{2, 7\} & \{4, 7\} & \{1, 5\} & \{3, 5\} & \{1, 2\} & \{6\} & \{3, 4\} \\
				\hline
				7	& \{2, 6\} & \{4, 6\} & \{2, 3\} & \{7\} & \{4, 5\} & \{1, 3\} & \{1, 5\} \\
				\hline
			\end{tabular}
		\end{center}
	\end{table}
	\newpage
	\par
	\noindent				
	{\bf Example 3.3:} (Example of polyloop)
	
	Let $G = \{1, 2, 3, 4, 5, 6\}$ and define hyperoperations $\odot,\nwarrow,\nearrow: G \times G \rightarrow \mathcal{P}^{*}(G)$ as in Table \ref{exampolyloop3}, Table \ref{exampolyloop4}, Table \ref{exampolyloop5}. Then, $(G, \odot, \nwarrow, \nearrow, 1)$, $(G, \nwarrow, \odot, \nwarrow,1)$ and $(G, \nearrow, \nearrow ,\odot,1)$ are polyloops.
	
	\begin{table}[!h]
		\begin{center}\caption{Multiplication table for  polyloop $(G, \odot, \nwarrow, \nearrow, 1)$}\label{exampolyloop3}
			\begin{tabular}{|c|c|c|c|c|c|c|}
				\hline
				$\odot$	& 1 & 2 & 3 & 4 & 5 & 6  \\
				\hline
				1	& \{1\} & \{2\} & \{3\} & \{4\} & \{5\} & \{6\}  \\
				\hline
				2	& \{2\} & \{1\} & \{4, 6\} & \{5, 6\} & \{3, 4\} & \{3, 5\}  \\
				\hline
				3	& \{3\} & \{4, 6\} & \{1\} & \{2, 5\} & \{2, 6\} & \{4, 5\}  \\
				\hline
				4	& \{4\} & \{5, 6\} & \{2, 5\} & \{1\} & \{3, 6\} & \{2, 3\}  \\
				\hline
				5	& \{5\} & \{3, 4\} & \{2, 6\} & \{3, 6\} & \{1\} & \{2, 4\}  \\
				\hline
				6	& \{6\} & \{3, 5\} & \{4, 5\} & \{2, 3\} & \{2, 4\} & \{1\}  \\
				\hline
			\end{tabular}
		\end{center}
	\end{table}
	\newpage
	\begin{table}[!h]
		\begin{center}\caption{Multiplication table for  polyloop $(G, \nwarrow, \odot, \nwarrow,1)$}\label{exampolyloop4}
			\begin{tabular}{|c|c|c|c|c|c|c|}
				\hline
				$\nwarrow$	& 1 & 2 & 3 & 4 & 5 & 6  \\
				\hline
				1	& \{1\} & \{2\} & \{3\} & \{4\} & \{5\} & \{6\}  \\
				\hline
				2	& \{2\} & \{1\} & \{5, 6\} & \{3, 5\} & \{4, 6\} & \{3, 4\}  \\
				\hline
				3	& \{3\} & \{4, 5\} & \{1\} & \{2, 6\} & \{4, 6\} & \{2, 5\}  \\
				\hline
				4	& \{4\} & \{3, 6\} & \{5, 6\} & \{1\} & \{2, 3\} & \{2, 5\}  \\
				\hline
				5	& \{5\} & \{3, 6\} & \{2, 4\} & \{2, 6\} & \{1\} & \{3, 4\}  \\
				\hline
				6	& \{6\} & \{4, 5\} & \{2, 4\} & \{3, 5\} & \{3, 4\} & \{1\}  \\
				\hline
				
			\end{tabular}
		\end{center}
	\end{table}
	
	\begin{table}[!h]
		\begin{center}\caption{Multiplication table for  polyloop $(G, \nearrow, \nearrow ,\odot,1)$ }\label{exampolyloop5}
			\begin{tabular}{|c|c|c|c|c|c|c|}
				\hline
				$\nearrow$	& 1 & 2 & 3 & 4 & 5 & 6  \\
				\hline
				1	& \{1\} & \{2\} & \{3\} & \{4\} & \{5\} & \{6\} \\
				\hline
				2	& \{2\} & \{1\} & \{4, 5\} & \{3, 6\} & \{3, 6\} & \{4, 5\}  \\
				\hline
				3	& \{3\} & \{5, 6\} & \{1\} & \{5, 6\} & \{2, 4\} & \{2, 4\}  \\
				\hline
				4	& \{4\} & \{3, 5\} & \{2, 6\} & \{1\} & \{2, 6\} & \{3, 5\}  \\
				\hline
				5	& \{5\} & \{4, 6\} & \{4, 6\} & \{2, 3\} & \{1\} & \{2, 3\}  \\
				\hline
				6	& \{6\} & \{3, 4\} & \{2, 5\} & \{2, 5\} & \{3, 4\} & \{1\}  \\
				\hline
				
			\end{tabular}
		\end{center}
	\end{table}
\subsection{Nuclei of Polygroupoid}
\paragraph{}
\vskip 0.3 cm	
We now introduce the notion of nuclei of polygroupoid as follows.
	\par
	\noindent				
	{\bf Definition 3.1:} (First Nuclei of Polygroupoid)
Let $(P, \cdot)$ be a polygroupoid. $N^1_{\lambda}(P, \cdot) \left(N^1_{\mu}(P, \cdot),N^1_{\rho}(P, \cdot)\right)$ will be called the first left (middle, right) nucleus of $(P, \cdot)$ and defined respectively as follows:

\begin{center}
	$N^1_{\lambda}(P, \cdot) = \{a \in P ~\mid~ a \cdot (x \cdot y) = (a \cdot x) \cdot y ~\forall\; x, y \in P\}$
	
	$N^1_{\mu}(P, \cdot) = \{a \in P ~\mid~ (x \cdot a) \cdot y = x \cdot (a \cdot y) ~\forall\; x, y \in P\}$
	
	$N^1_{\rho}(P, \cdot) = \{a \in P ~\mid~ (x \cdot y) \cdot a = x \cdot (y \cdot a) ~\forall\; x, y \in P\}$
\end{center}
The first nucleus of $(P, \cdot)$ will be given by

$N^1(P, \cdot) = N^1_{\lambda}(P, \cdot) \cap N^1_{\mu}(P, \cdot) \cap N^1_{\rho}(P, \cdot)$.	
	
\par
\noindent				
{\bf Definition 3.2:} (Second Nuclei of Polygroupoid)	
Let $(P, \cdot)$ be a polygroupoid. $N^2_{\lambda}(P, \cdot) \left(N^2_{\mu}(P, \cdot),N^2_{\rho}(P, \cdot)\right)$ will be called the second left (middle, right) nucleus of $(P, \cdot)$ defined respectively as follows:

\begin{center}
	$N^2_{\lambda}(P, \cdot) = \{a \in P ~\mid~ a \cdot (X \cdot y) = (a \cdot X) \cdot y ~\forall\; X\subseteq P,\; y \in P\}$
	
	$N^2_{\mu}(P, \cdot) = \{a \in P ~\mid~ (X \cdot a) \cdot y = X \cdot (a \cdot y) ~\forall\; X \subseteq P,\; y \in P\}$
	
	$N^2_{\rho}(P, \cdot) = \{a \in P ~\mid~ (X \cdot y) \cdot a = X \cdot (y \cdot a) ~\forall\; X \subseteq P,\; y \in P\}$
\end{center}
The second nucleus of $(P, \cdot)$ will be given by

$N^2(P, \cdot) = N^2_{\lambda}(P, \cdot) \cap N^2_{\mu}(P, \cdot) \cap N^2_{\rho}(P, \cdot)$.

\par
\noindent				
{\bf Definition 3.3:} (Third Nuclei of Polygroupoid)
Let $(P, \cdot)$ be polygroupoid. $N^3_{\lambda}(P, \cdot) \left(N^3_{\mu}(P, \cdot),N^3_{\rho}(P, \cdot)\right)$ will be called the third left (middle, right) nucleus of $(P, \cdot)$ defined respectively as follows:

\begin{center}
	$N^3_{\lambda}(P, \cdot) = \{a \in P ~\mid~ a \cdot (x \cdot Y) = (a \cdot x) \cdot Y ~\forall\; Y \subseteq P,\; x \in P\}$
	
	$N^3_{\mu}(P, \cdot) = \{a \in P ~\mid~ (x \cdot a) \cdot Y = x \cdot (a \cdot Y) ~\forall \;Y \subseteq P,\; x \in P\}$
	
	$N^3_{\rho}(P, \cdot) = \{a \in P ~\mid~ (x \cdot Y) \cdot a = x \cdot (Y \cdot a) ~\forall\; Y \subseteq P,\; x \in P\}$
\end{center}
The third nucleus of polygroupoid $(P, \cdot)$ will be given by

$N^3(P, \cdot) = N^3_{\lambda}(P, \cdot) \cap N^3_{\mu}(P, \cdot) \cap N^3_{\rho}(P, \cdot)$.

\par
\noindent				
{\bf Definition 3.4:} (Fourth Nuclei of Polygroupoid)
Let $(P, \cdot)$ be a polygroupoid. $N^4_{\lambda}(P, \cdot) \left(N^4_{\mu}(P, \cdot),N^4_{\rho}(P, \cdot)\right)$ will be called the fourth left (middle, right) nucleus of  $(P, \cdot)$ defined respectively as follows:

\begin{center}
	$N^4_{\lambda}(P, \cdot) = \{a \in P ~\mid~ a \cdot (X \cdot Y) = (a \cdot X) \cdot Y ~\forall \;X, Y \subseteq P\}$
	
	$N^4_{\mu}(P, \cdot) = \{a \in P ~\mid~ (X \cdot a) \cdot Y = X \cdot (a \cdot Y) ~\forall\; X, Y \subseteq P\}$
	
	$N^4_{\rho}(P, \cdot) = \{a \in P ~\mid~ (X \cdot Y) \cdot a = X \cdot (Y \cdot a) ~\forall\; X, Y \subseteq P\}$
\end{center}
The fourth nucleus of polygroupoid $(P, \cdot)$ will be given by

$N^4(P, \cdot) = N^4_{\lambda}(P, \cdot) \cap N^4_{\mu}(P, \cdot) \cap N^4_{\rho}(P, \cdot)$.

The relationship among the first, second, third and fourth the nuclei of polygroupoid are established as follows.	
	\par
	\noindent
	{\bf Theorem 3.1:} Let $(P, \cdot)$ be a polygroupoid and
	
	$N^1_{\lambda}(P, \cdot) = \{a \in P ~\mid~ a \cdot (x \cdot y) = (a \cdot x) \cdot y ~\forall x, y \in P\}$;
	
	$N^2_{\lambda}(P, \cdot) = \{a \in P ~\mid~ a \cdot (X \cdot y) = (a \cdot X) \cdot y ~\forall X \subseteq P, y \in P\}$;
	
	$N^3_{\lambda}(P, \cdot) = \{a \in P ~\mid~ a \cdot (x \cdot Y) = (a \cdot x) \cdot Y ~\forall Y \subseteq P, x \in P\}$;
	
	$N^4_{\lambda}(P, \cdot) = \{a \in P ~\mid~ a \cdot (X \cdot Y) = (a \cdot X) \cdot Y ~\forall X, Y \subseteq P\}$.
	
	Then, the following are true:
	\begin{multicols}{2}
		\begin{description}
			\item[(a)] $N^1_{\lambda}(P, \cdot) \subseteq N^2_{\lambda}(P, \cdot)$.
			\item[(b)] $N^1_{\lambda}(P, \cdot) \subseteq N^3_{\lambda}(P, \cdot)$.
			\item[(c)] $N^1_{\lambda}(P, \cdot) \subseteq N^4_{\lambda}(P, \cdot)$.
			\item[(d)] $N^2_{\lambda}(P, \cdot) \subseteq N^4_{\lambda}(P, \cdot)$.
			\item[(e)] $N^3_{\lambda}(P, \cdot) \subseteq N^4_{\lambda}(P, \cdot)$.
		\end{description}
	\end{multicols}
	\par
	\noindent
	{\bf Proof:} 	\begin{description}
		\item[(a)] Suppose $a \in N^1_{\lambda}(P, \cdot)$, then this implies $a \in P$ such that 
		
		$a \cdot (x \cdot y) = (a \cdot x) \cdot y ~ \forall~ x, y \in P$.
		
		Now,	
		\begin{equation}\label{eq1}
			\displaystyle\bigcup_{z \in x \cdot y}(a \cdot z) = \displaystyle\bigcup_{w \in a \cdot x}(w \cdot y)
		\end{equation}
		Thus, $z \in x \cdot y~\forall~x, y \in P$ implies $z \in \displaystyle\bigcup_{x \in X}(x \cdot y)$ and also, $w \in a \cdot x~\forall a \in P$ implies $w \in \displaystyle\bigcup_{x \in X}(a \cdot x)$. So, following \eqref{eq1}; 
		
		$\displaystyle\underset{\displaystyle z \in {\tiny \underset{x \in X}{\cup}}(x \cdot y)}{\bigcup (a \cdot z)} = \displaystyle\underset{\displaystyle w \in {\tiny \underset{x \in X}{\cup}}(a \cdot x)}{\bigcup (w \cdot y)}$
		
		$\Rightarrow~ \{a\} \cdot \displaystyle\bigcup_{x \in X}(x \cdot y) = \displaystyle\bigcup_{x \in X}(a \cdot x) \cdot \{y\}$
		
		$\Rightarrow~ a \cdot (X \cdot y) = (a \cdot X) \cdot y~ \forall ~ X \subseteq P, y \in P$
		
		$\Rightarrow~ a \in N^2_{\lambda}(P, \cdot)$.
		
		Hence, $N^1_{\lambda}(P, \cdot) \subseteq N^2_{\lambda}(P, \cdot)$.
		\item[(b)] Suppose $a \in N^1_{\lambda}(P, \cdot)$, then this implies $a \in P$ such that 
		
		$a \cdot (x \cdot y) = (a \cdot x) \cdot y ~ \forall~ x, y \in P$. Now,
		\begin{equation}\label{eq2}
			\displaystyle\bigcup_{z \in x \cdot y}(a \cdot z) = \displaystyle\bigcup_{w \in a \cdot x}(w \cdot y)
		\end{equation}
		So, $z \in x \cdot y~\forall~x, y \in P$ implies $z \in \displaystyle\bigcup_{y \in Y}(x \cdot y)$. Following \eqref{eq2}; 
		
		$\displaystyle\underset{\displaystyle z \in {\tiny \underset{y \in Y}{\cup}}(x \cdot y)}{\bigcup (a \cdot z)} = \displaystyle\bigcup_{y \in Y, w \in(a \cdot x)}(w \cdot y)$

		$\Rightarrow~ \{a\} \cdot \displaystyle\bigcup_{y \in Y}(x \cdot y) = (a \cdot x) \cdot Y$

		$\Rightarrow~ a \cdot (x \cdot Y) = (a \cdot x) \cdot Y~ \forall ~ Y \subseteq P, x \in P$

		$\Rightarrow~ a \in N^3_{\lambda}(P, \cdot)$. 

		Hence, $N^1_{\lambda}(P, \cdot) \subseteq N^3_{\lambda}(P, \cdot)$.		

		\item[(c)] Suppose $a \in N^1_{\lambda}(P, \cdot)$, this implies $a \in P$ such that 
		
		$a \cdot (x \cdot y) = (a \cdot x) \cdot y ~ \forall~ x, y \in P$. Now, 
		\begin{equation}\label{eq3}
			\displaystyle\bigcup_{z \in x \cdot y}(a \cdot z) = \displaystyle\bigcup_{w \in a \cdot x}(w \cdot y)
		\end{equation}
		Then, $z \in x \cdot y~\forall~x, y \in P$ implies $\displaystyle z \in \bigcup_{x \in X, y \in Y}(x \cdot y)$. If $w \in a \cdot x~\forall~ a \in P, x \in P$, this implies $w \in \displaystyle\bigcup_{x \in X}(a \cdot x)$. Following \eqref{eq3}; 
		
		$\displaystyle\underset{\displaystyle z \in {\tiny \underset{x \in X, y \in Y}{\cup}}(x \cdot y)}{\bigcup (a \cdot z)} = \displaystyle\underset{\displaystyle y \in Y, w \in {\tiny \underset{x \in X}{\cup}}(a \cdot y)}{\bigcup (w \cdot y)}$
		
		$\Rightarrow~ \{a\} \cdot \displaystyle\bigcup_{x \in X, y \in Y}(x \cdot y) = \displaystyle\bigcup_{x \in X}(a \cdot x) \cdot Y$
		
		$\Rightarrow~ a \cdot (X \cdot Y) = (a \cdot X) \cdot Y~ \forall ~ X, Y \subseteq P$
		
		$\Rightarrow~ a \in N^4_{\lambda}(P, \cdot)$
		
		Hence, $N^1_{\lambda} \subseteq N^4_{\lambda}(P, \cdot)$.		
		\item[(d)] Suppose $a \in N^2_{\lambda}(P, \cdot)$, then this implies $a \in P$ such that 
		
		$a \cdot (X \cdot y) = (a \cdot X) \cdot y ~ \forall~ x \subseteq P, y \in P$. 
		
		Now,
		\begin{equation}\label{eq4}
			\displaystyle\underset{\displaystyle z \in {\tiny \underset{x \in X}{\cup}}(x \cdot y)}{\bigcup (a \cdot z)} = \displaystyle\underset{\displaystyle w \in {\tiny \underset{x \in X}{\cup}}(a \cdot x)}{\bigcup (w \cdot y)}
		\end{equation}
		Then, if $z \in \displaystyle\bigcup_{x \in X}(x \cdot y)~\forall~ y \in P$, this implies $z \in \displaystyle\bigcup_{x \in X, y \in Y}(x \cdot y)$. Following \eqref{eq4}; 
		
		$\displaystyle\underset{\displaystyle z \in {\tiny \underset{x \in X, y \in Y}{\cup}}(x \cdot y)}{\bigcup (a \cdot z)} = \displaystyle\underset{\displaystyle y \in Y, z \in {\tiny \underset{x \in X}{\cup}}(a \cdot x)}{\bigcup (w \cdot y)}$
		
		$\Rightarrow~ \{a\} \cdot \displaystyle\bigcup_{x \in X, y \in Y}(x \cdot y) = \displaystyle\bigcup_{x \in X}(a \cdot x) \cdot Y$
		
		$\Rightarrow~ a \cdot (X \cdot Y) = (a \cdot X) \cdot Y~ \forall ~ X, Y \subseteq P$
		
		$\Rightarrow~ a \in N^4_{\lambda}(P, \cdot)$.
		
		Hence, $N^2_{\lambda}(P, \cdot) \subseteq N^4_{\lambda}(P, \cdot)$.
		\item[(e)] Suppose $a \in N^3_{\lambda}(P, \cdot)$, then this implies $a \in P$ such that 
		
		$a \cdot (x \cdot Y) = (a \cdot x) \cdot Y ~ \forall~ Y \subseteq P, x \in P$. 
		
		Now,
		\begin{equation}\label{eq5}
			\displaystyle\underset{\displaystyle z \in {\tiny \underset{y \in Y}{\cup}}(x \cdot y)}{\bigcup (a \cdot z)} = \displaystyle\bigcup_{w \in (a \cdot x), y \in Y}(w \cdot y)
		\end{equation}
		Then, if $z \in \displaystyle\bigcup_{y \in Y}(x \cdot y)~\forall~ x \in P$, this implies $z \in \displaystyle\bigcup_{x \in X, y \in Y}(x \cdot y)$. If $w \in (a \cdot x)$, then this implies $w \in \displaystyle\bigcup_{x \in X}(a \cdot x)$. Following \eqref{eq5}; 
		
		$\displaystyle\underset{\displaystyle z \in {\tiny \underset{x \in X, y \in Y}{\cup}}(x \cdot y)}{\bigcup (a \cdot z)} = \displaystyle\underset{\displaystyle y \in Y, w \in {\tiny \underset{x \in X}{\cup}}(a \cdot x)}{\bigcup (w \cdot y)}$
		
		$\Rightarrow~ \{a\} \cdot \displaystyle\bigcup_{x \in X, y \in Y}(x \cdot y) = \displaystyle\bigcup_{x \in X}(a \cdot x) \cdot Y$
		
		$\Rightarrow~ a \cdot (X \cdot Y) = (a \cdot X) \cdot Y~ \forall ~ X, Y \subseteq P$
		
		$\Rightarrow~ a \in N^4_{\lambda}(P, \cdot)$.
		
		Hence, $N^3_{\lambda}(P, \cdot) \subseteq N^4_{\lambda}(P, \cdot)$.
	\end{description}
	
		\par
	\noindent
	{\bf Remark 3.1:} In a groupoid (quasigroup, loop), the left nucleus coincides with the first left nucleus of  polygroupoid (polyquasigroup, polyloop) in Theorem $3.1$ while second, third and fourth left nuclei of  polygroupoid (polyquasigroup, polyloop) will only coincide in a groupoid (quasigroup, loop) if the arbitrary elements in the power set in the axioms are replaced with arbitrary singleton sets.
	
	\par 
	\noindent
	{\bf Remark 3.2:} The first left nucleus polygroupoid (polyquasigroup, polyloop) in Theorem $3.1$ is contained in second and third left nuclei of polygroupoid (polyquasigroup, polyloop) in Theorem $3.1$. This shows that second and third left nuclei of a polygroupoid (polyquasigroup, polyloop) generalize the first left nucleus of polygroupoid (polyquasigroup, polyloop). Furthermore, the fourth left nucleus of polyloop generalizes first, second and third left nuclei of polygroupoid (polyquasigroup, polyloop).

	\par 
	\noindent
	{\bf Theorem 3.2:} Let $(P, \cdot)$ be a polygroupoid and 
	
	$N^1_{\mu}(P, \cdot) = \{a \in P ~\mid~ (x \cdot a) \cdot y = x \cdot (a \cdot y) ~\forall~ x, y \in P\}$;
	
	$N^2_{\mu}(P, \cdot) = \{a \in P ~\mid~ (X \cdot a) \cdot y = X \cdot (a \cdot y) ~\forall~ X \subseteq P, y \in P\}$;
	
	$N^3_{\mu}(P, \cdot) = \{a \in P ~\mid~ (x \cdot a) \cdot Y = x \cdot (a \cdot Y) ~\forall~ x \in P, Y \subseteq P\}$;
	
	$N^4_{\mu}(P, \cdot) = \{a \in P ~\mid~ (X \cdot a) \cdot Y = X \cdot (a \cdot Y) ~\forall~ X, Y \subseteq P\}$.
	
	Then, the following are true:

	\begin{multicols}{2}
		\begin{description}
					\item[(a)] $N^1_{\mu}(P, \cdot) \subseteq N^2_{\mu}(P, \cdot)$.
			\item[(b)] $N^1_{\mu}(P, \cdot) \subseteq N^3_{\mu}(P, \cdot)$.
			\item[(c)] $N^1_{\mu}(P, \cdot) \subseteq N^4_{\mu}(P, \cdot)$.
			\item[(d)] $N^2_{\mu}(P, \cdot) \subseteq N^4_{\mu}(P, \cdot)$.
			\item[(e)] $N^3_{\mu}(P, \cdot) \subseteq N^4_{\mu}(P, \cdot)$.
		\end{description}
	\end{multicols}
	
\par
\noindent
{\bf Proof:} \begin{description}
	\item[(a)] Suppose $a \in N^1_{\mu}(P, \cdot)$,  then $(x \cdot a) \cdot y = x \cdot (a \cdot y) ~ \forall~ x, y \in P$. Thus,
	\begin{equation}\label{eq6}
		\displaystyle\bigcup_{w \in x \cdot a}(w \cdot y) = \displaystyle\bigcup_{z \in a \cdot y}(x \cdot z)
	\end{equation}
	If $w \in x \cdot a~\forall~ x \in P$, this implies $w \in \displaystyle\bigcup_{x \in X}(x \cdot a)$. Following \eqref{eq6}; 
	
	$\displaystyle\underset{\displaystyle w \in {\tiny \underset{x \in X}{\cup}}(x \cdot a)}{\bigcup (w \cdot y)} = \displaystyle\bigcup_{z \in a \cdot y, x \in X}(x \cdot z)$
	
	$\Rightarrow \left(\displaystyle\bigcup_{x \in X}(x \cdot a)\right) \cdot \{y\} = X \cdot (a \cdot y)$
	
	$\Rightarrow~ (X \cdot a) \cdot y = X \cdot (a \cdot y)$
	
	$\Rightarrow~ a \in N^2_{\mu}(P, \cdot)$.
	
	Hence, $N^1_{\mu} \subseteq N^2_{\mu}(P, \cdot)$.
	\item[(b)] Suppose $a \in N^1_{\mu}(P, \cdot)$, then $(x \cdot a) \cdot y = x \cdot (a \cdot y) ~ \forall~ x, y \in P$. Now,
	\begin{equation}\label{eq7}
		\Rightarrow~ \displaystyle\bigcup_{w \in x \cdot a}(w \cdot y) = \displaystyle\bigcup_{z \in a \cdot y}(x \cdot z)
	\end{equation}
	If $z \in a \cdot y~\forall~ y \in P$, then this implies $w \in \displaystyle\bigcup_{y \in Y}(a \cdot y)$. Following \eqref{eq7}; 
	
	$\displaystyle\bigcup_{w \in (x \cdot a), y \in Y}(w \cdot y) = \displaystyle\underset{\displaystyle z \in {\tiny \underset{y \in Y}{\cup}}(a \cdot y)}{\bigcup (x \cdot z)}$
	
	$\Rightarrow ~(x \cdot a) \cdot Y = \{x\} \cdot \displaystyle\bigcup_{y \in Y}(a \cdot y)$
	
	$\Rightarrow~ (x \cdot a) \cdot Y = x \cdot (a \cdot Y) ~ \forall~ x \in P, Y \subseteq P$
	
	$\Rightarrow~ a \in N^3_{\mu}(P, \cdot)$. 
	
	Hence,$N^1_{\mu}(P, \cdot) \subseteq N^3_{\mu}(P, \cdot)$.
	\item[(c)] Suppose $a \in N^1_{\mu}(P, \cdot)$, then $(x \cdot a) \cdot y = x \cdot (a \cdot y) ~ \forall~ x, y \in P$. Now, 
	\begin{equation}\label{eq8}
		\displaystyle\bigcup_{w \in x \cdot a}(w \cdot y) = \displaystyle\bigcup_{z \in a \cdot y}(x \cdot z)
	\end{equation}
	If $w \in x \cdot a$, then this implies $w \in \displaystyle\bigcup_{x \in X}(x \cdot a)$. 
	If $z \in a \cdot y~\forall~ y \in P$, then this implies $w \in \displaystyle\bigcup_{y \in Y}(a \cdot y)$. Following \eqref{eq8}; 
	
	$\displaystyle\bigcup_{w \in (x \cdot a), y \in Y}(w \cdot y) = \displaystyle\underset{\displaystyle z \in {\tiny \underset{y \in Y}{\cup}}(a \cdot y), x \in X}{\bigcup (x \cdot z)}$
	
	$\Rightarrow ~\displaystyle\bigcup_{x \in X}(x \cdot a) \cdot Y = X \cdot \displaystyle\bigcup_{y \in Y}(a \cdot y)$
	
	$\Rightarrow~ (X \cdot a) \cdot Y = X \cdot (a \cdot Y) ~ \forall~ X, Y \subseteq P$
	
	$\Rightarrow~ a \in N^4_{\mu}(P, \cdot)$. Hence, $N^1_{\mu}(P, \cdot) \subseteq N^4_{\mu}(P, \cdot)$.		
	\item[(d)] Suppose $a \in N^2_{\mu}(P, \cdot)$, then this  implies $a \in P$ such that 
	
	$(X \cdot a) \cdot y = X \cdot (a \cdot y) ~ \forall~ X \subseteq P, y \in P$. Now,
	\begin{equation}\label{eq9}
		\displaystyle\underset{\displaystyle w \in {\tiny \underset{x\in X}{\cup}}(x \cdot a)}{\bigcup (w \cdot y)} = \displaystyle\bigcup_{z \in a \cdot y, x \in X}(x \cdot z)
	\end{equation}
	Following \eqref{eq9}; 
	
	$\displaystyle\underset{\displaystyle w \in {\tiny \underset{x \in X}{\cup}}(a \cdot y), y \in Y}{\bigcup (w \cdot y)} = \displaystyle\underset{\displaystyle z \in {\tiny \underset{y \in Y}{\cup}}(a \cdot y), x \in X}{\bigcup (x \cdot z)}$
	
	$\Rightarrow ~\displaystyle\bigcup_{x \in X}(x \cdot a) \cdot Y = X \cdot \displaystyle\bigcup_{y \in Y}(a \cdot y)$
	
	$\Rightarrow~ (X \cdot a) \cdot Y = X \cdot (a \cdot Y) ~ \forall~ X, Y \subseteq P$
	
	$\Rightarrow~ a \in N^4_{\mu}(P, \cdot)$.
	
	Hence, $N^2_{\mu}(P, \cdot) \subseteq N^4_{\mu}(P, \cdot)$.		
	
	\item[(e)] Suppose $a \in N^3_{\mu}(P, \cdot)$, then this implies $a \in P$ such that 
	
	$(x \cdot a) \cdot Y = x \cdot (a \cdot Y) ~ \forall~ Y \subseteq P, x \in P$. Now,
	\begin{equation}\label{eq10}
		\displaystyle\bigcup_{w \in (x \cdot a), y \in Y}(w \cdot y) = \displaystyle\underset{\displaystyle z \in {\tiny \underset{y \in Y}{\cup}}(a \cdot y)}{\bigcup (x \cdot z)}
	\end{equation}
	If $w \in (x \cdot a)~\forall~ x \in P$, then this implies $w \in \displaystyle\bigcup_{x \in X}(x \cdot a)$. Following \eqref{eq10}; 
	
	$\displaystyle\underset{\displaystyle w \in {\tiny \underset{x \in X}{\cup}}(a \cdot y), y \in Y}{\bigcup (x \cdot z)} = \displaystyle\underset{\displaystyle z \in {\tiny \underset{y \in Y}{\cup}}(a \cdot y), x \in X}{\bigcup (x \cdot z)}$
	
	$\Rightarrow ~\displaystyle\bigcup_{x \in X}(x \cdot a) \cdot Y = X \cdot \displaystyle\bigcup_{y \in Y}(a \cdot y)$
	
	$\Rightarrow~ (X \cdot a) \cdot Y = X \cdot (a \cdot Y) ~ \forall~ X, Y \subseteq P$
	
	$\Rightarrow~ a \in N^4_{\mu}(P, \cdot)$.
	
	Hence, $N^3_{\mu}(P, \cdot) \subseteq N^4_{\mu}(P, \cdot)$.
\end{description}	
	
\par 
\noindent
{\bf Remark 3.3:}	In a groupoid (quasigroup, loop), the  middle nucleus coincides with the first middle nucleus of  polygroupoid (polyquasigroup, polyloop) in Theorem $3.2$ while second, third and fourth middle nuclei of  polygroupoid (polyquasigroup, polyloop) will only coincide in a groupoid (quasigroup, loop) if the arbitrary elements in the power set in the axioms are replaced with arbitrary singleton sets.
	
\par 
\noindent
{\bf Remark 3.4:} The first  middle nucleus polygroupoid (polyquasigroup, polyloop) in Theorem $3.2$ is contained in second and third  middle nuclei of polygroupoid (polyquasigroup, polyloop) in Theorem $3.2$. This shows that second and third middle nuclei of a polygroupoid (polyquasigroup, polyloop) generalize the first  middle nucleus of polygroupoid (polyquasigroup, polyloop). Furthermore, the fourth  middle nucleus of polyloop generalizes first, second and third  middle nuclei of polygroupoid (polyquasigroup, polyloop).

	\par 
\noindent
{\bf Theorem 3.3:} \label{myth3} Let $(P, \cdot)$ be a polygroupoid and 

$N^1_{\rho}(P, \cdot) = \{a \in P ~\mid~ (x \cdot y) \cdot a =  \cdot (y \cdot a) ~\forall~ x, y \in P\}$;

$N^2_{\rho}(P, \cdot) = \{a \in P ~\mid~ (X \cdot y) \cdot a = X \cdot (y \cdot a) ~\forall~ X \subseteq P, y \in P\}$;

$N^3_{\rho}(P, \cdot) = \{a \in P ~\mid~ (x \cdot Y) \cdot a = x \cdot (Y \cdot a) ~\forall~ x \in P, Y \subseteq P\}$;

$N^4_{\rho}(P, \cdot) = \{a \in P ~\mid~ (X \cdot Y) \cdot a = X \cdot (Y \cdot a) ~\forall~ X, Y \subseteq P\}$.

Then the following are true:
\begin{multicols}{2}
	\begin{description}
		\item[(a)] $N^1_{\rho}(P, \cdot) \subseteq N^2_{\rho}(P, \cdot)$.
		\item[(b)] $N^1_{\rho}(P, \cdot) \subseteq N^3_{\rho}(P, \cdot)$.
		\item[(c)] $N^1_{\rho}(P, \cdot) \subseteq N^4_{\rho}(P, \cdot)$.
		\item[(d)] $N^2_{\rho}(P, \cdot) \subseteq N^4_{\rho}(P, \cdot)$.
		\item[(e)] $N^3_{\rho}(P, \cdot) \subseteq N^4_{\rho}(P, \cdot)$.	 
	\end{description}
\end{multicols}

	\par 
\noindent
{\bf Proof:} \begin{description}
	\item[(a)] Suppose $a \in N^1_{\rho}(P, \cdot)$, then this  implies $a \in P$ such that 
	\begin{center}
	$(x \cdot y) \cdot a = x \cdot (y \cdot a) ~ \forall~ x, y \in P$.
	\end{center} 
	 Now,
	
	\begin{equation}\label{eq11}
		\displaystyle\bigcup_{b \in x \cdot y}(b \cdot a) = \displaystyle\bigcup_{w \in y \cdot a}(x \cdot w)
	\end{equation}
	If $b \in x \cdot y~ \forall ~x, y \in P$, then this implies $b \in \displaystyle\bigcup_{x \in X}(x \cdot y)$. 	
	Following \eqref{eq11}; $\displaystyle\underset{\displaystyle b \in {\tiny \underset{x \in X}{\cup}}(x \cdot y)}{\bigcup (b \cdot a)} = \displaystyle\bigcup_{w \in y \cdot a, x \in X}(x \cdot w)$
	
	$\Rightarrow~ (\displaystyle\bigcup_{x \in X}(x \cdot y)) \cdot \{a\} = X \cdot (y \cdot a)$
	
	$\Rightarrow~ (X \cdot y) \cdot a = X \cdot (y \cdot a) ~\forall~ y \in P, X \subseteq P$
	
	$\Rightarrow ~a \in N^2_{\rho}(P, \cdot)$.
	
	Hence, $N^1_{\rho}(P, \cdot) \subseteq N^2_{\rho}(P, \cdot)$.
	\item[(b)] Suppose $a \in N^1_{\rho}(P, \cdot)$, then implies $a \in P$ such that 
	
	$(x \cdot Y) \cdot a = x \cdot (Y \cdot a) ~ \forall~ x \in P, Y \subseteq P$.
	\begin{equation}\label{eq12}
		\displaystyle\bigcup_{b \in x \cdot y}(b \cdot a) = \displaystyle\bigcup_{w \in y \cdot a}(x \cdot w)
	\end{equation}
	If $b \in x \cdot y~ \forall ~x, y \in P$, then this implies $b \in \displaystyle\bigcup_{y \in Y}(x \cdot y)$. 		
	If $w \in y \cdot a ~\forall~ y \in P$ implies $w \in \displaystyle\bigcup_{y \in Y}(y \cdot a)$. Following \eqref{eq12}; $\displaystyle\underset{\displaystyle b \in {\tiny \underset{y \in Y}{\cup}}(x \cdot y)}{\bigcup (b \cdot a)} = \displaystyle\underset{\displaystyle w \in {\tiny \underset{y \in Y}{\cup}}(y \cdot a)}{\bigcup (x \cdot w)}$
	
	$\Rightarrow~ (\displaystyle\bigcup_{y \in Y}(x \cdot y)) \cdot \{a\} = \{x\} \cdot \displaystyle\bigcup_{y \in Y}(y \cdot a)$
	
	$\Rightarrow~ (x \cdot Y) \cdot a = x \cdot (Y \cdot a) ~\forall~ x \in P, Y \subseteq P$
	
	$\Rightarrow ~a \in N^3_{\rho}(P, \cdot)$.
	
	Hence, $N^1_{\rho}(P, \cdot) \subseteq N^3_{\rho}(P, \cdot)$.		
	
	\item[(c)] Suppose $a \in N^1_{\rho}(P, \cdot)$, then this implies $a \in P$ such that 
	
	$(x \cdot y) \cdot a = x \cdot (y \cdot a) ~ \forall~ x, y \in P$. Now, 
	\begin{equation}\label{eq13}
		\displaystyle\bigcup_{b \in x \cdot y}(b \cdot a) = \displaystyle\bigcup_{w \in y \cdot a}(x \cdot w)
	\end{equation}
	If $b \in x \cdot y~ \forall ~x, y \in P$, then this implies $b \in \displaystyle\bigcup_{x \in X, y \in Y}(x \cdot y)$. 	
	If $w \in y \cdot a ~\forall~ y \in P$, then implies $w \in \displaystyle\bigcup_{y \in Y}(y \cdot a)$. 	
	Following \eqref{eq13}; $\displaystyle\underset{\displaystyle b \in {\tiny \underset{x \in X, y \in Y}{\cup}}(x \cdot y)}{\bigcup (b \cdot a)} = \displaystyle\underset{\displaystyle w \in {\tiny \underset{y \in Y}{\cup}}(x \cdot y), x \in X}{\bigcup (x \cdot w)}$
	
	$\Rightarrow~ \displaystyle\bigcup_{x \in X, y \in Y}(x \cdot y) \cdot \{a\} = X \cdot \displaystyle\bigcup_{y \in Y}(y \cdot a)$
	
	$\Rightarrow~ (X \cdot Y) \cdot a = X \cdot (Y \cdot a) ~\forall~ X, Y \subseteq P$
	
	$\Rightarrow ~a \in N^4_{\rho}(P, \cdot)$.
	
	Hence, $N^1_{\rho}(P, \cdot) \subseteq N^4_{\rho}(P, \cdot)$.		
	\item[(d)] Suppose $a \in N^2_{\rho}(P, \cdot)$, then this implies $a \in P$ such that 
	
	$(X \cdot y) \cdot a = X \cdot (y \cdot a) ~ \forall~ y \in P, X \subseteq P$. Now,
	\begin{equation}\label{eq14}
		\displaystyle\underset{\displaystyle b \in {\tiny \underset{x \in X}{\cup}}(x \cdot y)}{\bigcup (b \cdot a)} = \displaystyle\bigcup_{w \in y \cdot a, x \in X}(x \cdot w)
	\end{equation}
	If $b \in \displaystyle\bigcup_{x \in X}(x \cdot y)~ \forall ~ y \in P$, then this implies $b \in \displaystyle\bigcup_{x \in X, y \in Y}(x \cdot y)$. Following \eqref{eq14}; $\displaystyle\underset{\displaystyle b \in {\tiny \underset{x \in X, y \in Y}{\cup}}(x \cdot y)}{\bigcup (b \cdot a)} = \displaystyle\underset{\displaystyle w \in {\tiny \underset{y \in Y}{\cup}}(y \cdot a), x \in X}{\bigcup (x \cdot w)}$
	
	$\Rightarrow~ \displaystyle\bigcup_{x \in X, y \in Y}(x \cdot y) \cdot \{a\} = X \cdot \displaystyle\bigcup_{y \in Y}(y \cdot a)$
	
	$\Rightarrow~ (X \cdot Y) \cdot a = X \cdot (Y \cdot a) ~\forall~ X, Y \subseteq P$
	
	$\Rightarrow ~a \in N^4_{\rho}(P, \cdot)$.
	
	Hence, $N^2_{\rho}(P, \cdot) \subseteq N^4_{\rho}(P, \cdot)$.
	\item[(e)] Suppose $a \in N^3_{\rho}(P, \cdot)$, then this implies $a \in P$ such that 
	
	$(x \cdot Y) \cdot a = x \cdot (Y \cdot a) ~ \forall~ x \in P, Y \subseteq P$. Now, 
	\begin{equation}\label{eq15}
		\displaystyle\underset{\displaystyle b \in {\tiny \underset{x \in X}{\cup}}(x \cdot y)}{\bigcup (b \cdot a)} = \displaystyle\underset{\displaystyle w \in {\tiny \underset{y \in Y}{\cup}}(x \cdot y)}{\bigcup (x \cdot w)}
	\end{equation}
	If $b \in \displaystyle\bigcup_{y \in Y}(x \cdot y)~ \forall ~ y \in P$, then this implies $b \in \displaystyle\bigcup_{x \in X, y \in Y}(x \cdot y)$. Following \eqref{eq15}; $\displaystyle\underset{\displaystyle b \in {\tiny \underset{x \in X, y \in Y}{\cup}}(x \cdot y)}{\bigcup (b \cdot a)} = \displaystyle\underset{\displaystyle w \in {\tiny \underset{y \in Y}{\cup}}(x \cdot y), x \in X}{\bigcup (x \cdot w)}$
	
	$\Rightarrow~ \displaystyle\bigcup_{x \in X, y \in Y}(x \cdot y) \cdot \{a\} = X \cdot \displaystyle\bigcup_{y \in Y}(y \cdot a)$
	
	$\Rightarrow~ (X \cdot Y) \cdot a = X \cdot (Y \cdot a) ~\forall~ X, Y \subseteq P$
	
	$\Rightarrow ~a \in N^4_{\rho}(P, \cdot)$.
	
	Hence, $N^3_{\rho}(P, \cdot) \subseteq N^4_{\rho}(P, \cdot)$.
\end{description}

\par 
\noindent
{\bf Remark 3.5:} In a groupoid (quasigroup, loop), the  right nucleus coincides with the first right nucleus of  polygroupoid (polyquasigroup, polyloop) in Theorem \ref{myth3} while second, third and fourth right nuclei of  polygroupoid (polyquasigroup, polyloop) will only coincide in a groupoid (quasigroup, loop) if the arbitrary elements in the power set in the axioms are replaced with arbitrary singleton sets.

\par 
\noindent
{\bf Remark 3.6:} The first  right nucleus polygroupoid (polyquasigroup, polyloop) in Theorem \ref{myth3} is contained in second and third  right nuclei of polygroupoid (polyquasigroup, polyloop) in Theorem \ref{myth3}. This shows that second and third right nuclei of a polygroupoid (polyquasigroup, polyloop) generalize the first  right nucleus of polygroupoid (polyquasigroup, polyloop). Furthermore, the fourth  right nucleus of polyloop generalizes first, second and third right nuclei of polygroupoid (polyquasigroup, polyloop).

\par 
\noindent
{\bf Theorem 3.4:} Let $N^1(P, \cdot), N^2(P, \cdot), N^3(P, \cdot), N^4(P, \cdot)$ denote the first, second, third and fourth nucleus of polygroupoid $(P, \cdot)$ respectively such that

$N^1(P, \cdot) = N^1_{\lambda}(P, \cdot) \cap N^1_{\mu}(P, \cdot) \cap N^1_{\rho}(P, \cdot)$

$N^2(P, \cdot) = N^2_{\lambda}(P, \cdot) \cap N^2_{\mu}(P, \cdot) \cap N^2_{\rho}(P, \cdot)$

$N^3(P, \cdot) = N^3_{\lambda}(P, \cdot) \cap N^3_{\mu}(P, \cdot) \cap N^3_{\rho}(P, \cdot)$

$N^4(P, \cdot) = N^4_{\lambda}(P, \cdot) \cap N^4_{\mu}(P, \cdot) \cap N^4_{\rho}(P, \cdot)$

Then the following are true:
\begin{multicols}{2}
	\begin{description}
		\item[(i)] $N^1(P, \cdot) \subseteq N^2(P, \cdot)$.
		\item[(ii)] $N^1(P, \cdot) \subseteq N^3(P, \cdot)$.
		\item[(iii)] $N^1(P, \cdot) \subseteq N^4(P, \cdot)$.
		\item[(iv)] $N^2(P, \cdot) \subseteq N^4(P, \cdot)$.
		\item[(v)] $N^3(P, \cdot) \subseteq N^4(P, \cdot)$.
	\end{description}
\end{multicols}
\par 
\noindent
{\bf Proof:} Let $N^1(P, \cdot), N^2(P, \cdot), N^3(P, \cdot), N^4(P, \cdot)$ respectively denote the first, second, third and fourth nuclei of polygroupoid $(P, \cdot)$. 
\begin{description}
	\item[(i)] $N^1(P, \cdot) = N^1_{\lambda}(P, \cdot) \cap N^1_{\mu}(P, \cdot) \cap N^1_{\rho}(P, \cdot)$ and $N^2(P, \cdot) = N^2_{\lambda}(P, \cdot) \cap N^2_{\mu}(P, \cdot) \cap N^2_{\rho}(P, \cdot)$ are first and second nucleus of polygroupoid $(P, \cdot)$ and we are to show that $N(P, \cdot) \subseteq N^2(P, \cdot)$.
	
	Suppose $x \in N^1(P, \cdot)$, this implies that $x \in N^1_{\lambda}(P, \cdot) \cap N^1_{\mu}(P, \cdot) \cap N^1_{\rho}(P, \cdot)$
	
	$\Rightarrow~ x \in N^1_{\lambda}(P, \cdot), x \in N^1_{\mu}(P, \cdot), x \in N^1_{\rho}(P, \cdot)$.
	
	Since, $N^1_{\lambda}(P, \cdot) \subseteq N^2_{\lambda}(P, \cdot), N^1_{\mu}(P, \cdot) \subseteq N^2_{\mu}(P, \cdot), N^1_{\rho}(P, \cdot) \subseteq N^2_{\rho}(P, \cdot)$
	$\Rightarrow~ x \in N^2_{\lambda}(P, \cdot) \cap N^2_{\mu}(P, \cdot) \cap N^2_{\rho}(P, \cdot)$
	Hence, $N^1(P, \cdot) \subseteq N^2(P, \cdot)$.
	\item[(ii)] $N^1(P, \cdot) = N^1_{\lambda}(P, \cdot) \cap N^1_{\mu}(P, \cdot) \cap N^1_{\rho}(P, \cdot)$ and $N^3(P, \cdot) = N^3_{\lambda}(P, \cdot) \cap N^3_{\mu}(P, \cdot) \cap N^3_{\rho}(P, \cdot)$ are respectively the first and third nuclei of polygroupoid $(P, \cdot)$ and we are to show that to show that $N^1(P, \cdot) \subseteq N^3(P, \cdot)$.
	
	Suppose $x \in N^1(P, \cdot)$ implies that $x \in N^1_{\lambda}(P, \cdot) \cap N^1_{\mu}(P, \cdot) \cap N^1_{\rho}(P, \cdot)$
	$\Rightarrow~ x \in N^1_{\lambda}(P, \cdot), x \in N^1_{\mu}(P, \cdot), x \in N^1_{\rho}(P, \cdot)$.
	
	Since, $N^1_{\lambda}(P, \cdot) \subseteq N^3_{\lambda}(P, \cdot), N^1_{\mu}(P, \cdot) \subseteq N^3_{\mu}(P, \cdot), N^1_{\rho}(P, \cdot) \subseteq N^3_{\rho}(P, \cdot)$
	$\Rightarrow~ x \in N^3_{\lambda}(P, \cdot) \cap N^3_{\mu}(P, \cdot) \cap N^3_{\rho}(P, \cdot)$.
	Hence, $N^1(P, \cdot) \subseteq N^3(P, \cdot)$.
	
	\item[(iii)] $N^1(P, \cdot) = N^1_{\lambda}(P, \cdot) \cap N^1_{\mu}(P, \cdot) \cap N^1_{\rho}(P, \cdot)$ and $N^4(P, \cdot) = N^4_{\lambda}(P, \cdot) \cap N^4_{\mu}(P, \cdot) \cap N^4_{\rho}(P, \cdot)$ are respectively the first and fourth nuclei of polygroupoid $(P, \cdot)$ and we are to show that $N^1(P, \cdot) \subseteq N^4(P, \cdot)$.
	
	Suppose $x \in N^1(P, \cdot)$ implies that $x \in N^1_{\lambda}(P, \cdot) \cap N^1_{\mu}(P, \cdot) \cap N^1_{\rho}(P, \cdot)$
	$\Rightarrow~ x \in N^1_{\lambda}(P, \cdot), x \in N^1_{\mu}(P, \cdot), x \in N^1_{\rho}(P, \cdot)$.
	
	Since, $N^1_{\lambda}(P, \cdot) \subseteq N^4_{\lambda}(P, \cdot), N^1_{\mu}(P, \cdot) \subseteq N^4_{\mu}(P, \cdot), N^1_{\rho}(P, \cdot) \subseteq N^4_{\rho}(P, \cdot)$
	$\Rightarrow~ x \in N^4_{\lambda}(P, \cdot) \cap N^4_{\mu}(P, \cdot) \cap N^4_{\rho}(P, \cdot)$. Hence, $N^1(P, \cdot) \subseteq N^4(P, \cdot)$.
	
	\item[(iv)] $N^2(P, \cdot) = N^2_{\lambda}(P, \cdot) \cap N^2_{\mu}(P, \cdot) \cap N^2_{\rho}(P, \cdot)$ and $N^4(P, \cdot) = N^4_{\lambda}(P, \cdot) \cap N^4_{\mu}(P, \cdot) \cap N^4_{\rho}(P, \cdot)$ are second and fourth nucleus of polygroupoid $(P, \cdot)$ and we are to show that $N^2(P, \cdot) \subseteq N^4(P, \cdot)$.
	
	Suppose $x \in N(P, \cdot)$ implies that $x \in N^2_{\lambda}(P, \cdot) \cap N^2_{\mu}(P, \cdot) \cap N^2_{\rho}(P, \cdot)$
	$\Rightarrow~ x \in N^2_{\lambda}(P, \cdot), x \in N^2_{\mu}(P, \cdot), x \in N^2_{\rho}(P, \cdot)$.
	
	Since, $N^2_{\lambda}(P, \cdot) \subseteq N^4_{\lambda}(P, \cdot), N^2_{\mu}(P, \cdot) \subseteq N^4_{\mu}(P, \cdot), N^2_{\rho}(P, \cdot) \subseteq N^4_{\rho}(P, \cdot)$
	$\Rightarrow~ x \in N^4_{\lambda}(P, \cdot) \cap N^4_{\mu}(P, \cdot) \cap N^4_{\rho}(P, \cdot)$. Hence, $N^2(P, \cdot) \subseteq N^4(P, \cdot)$.
	
	\item[(v)] $N^3(P, \cdot) = N^3_{\lambda}(P, \cdot) \cap N^3_{\mu}(P, \cdot) \cap N^3_{\rho}(P, \cdot)$ and $N^4(P, \cdot) = N^4_{\lambda}(P, \cdot) \cap N^4_{\mu}(P, \cdot) \cap N^4_{\rho}(P, \cdot)$ are third and fourth nucleus of polygroupoid $(P, \cdot)$ and we are to show that $N^3(P, \cdot) \subseteq N^4(P, \cdot)$.
	
	Suppose $x \in N^3(P, \cdot)$ implies that $x \in N^3_{\lambda}(P, \cdot) \cap N^3_{\mu}(P, \cdot) \cap N^3_{\rho}(P, \cdot)$
	$\Rightarrow~ x \in N^3_{\lambda}(P, \cdot), x \in N^3_{\mu}(P, \cdot), x \in N^3_{\rho}(P, \cdot)$.
	
	Since, $N^3_{\lambda}(P, \cdot) \subseteq N^4_{\lambda}(P, \cdot), N^3_{\mu}(P, \cdot) \subseteq N^4_{\mu}(P, \cdot), N^3_{\rho}(P, \cdot) \subseteq N^4_{\rho}(P, \cdot)$
	$\Rightarrow~ x \in N^4_{\lambda}(P, \cdot) \cap N^4_{\mu}(P, \cdot) \cap N^4_{\rho}(P, \cdot)$. Hence, $N^3(P, \cdot) \subseteq N^4(P, \cdot)$.
\end{description}

\par 
\noindent
{\bf Remark 3.7:} In a groupoid (quasigroup, loop), the   nucleus coincides with the first nucleus of  polygroupoid (polyquasigroup, polyloop) in Theorem $3.4$ while second, third and fourth nuclei of  polygroupoid (polyquasigroup, polyloop) will only coincide in a groupoid (quasigroup, loop) if the arbitrary elements in the power set in the axioms are replaced with arbitrary singleton sets.

\par 
\noindent
{\bf Remark 3.8:} The first  nucleus polygroupoid (polyquasigroup, polyloop) in Theorem $3.4$ is contained in second and third   nuclei of polygroupoid (polyquasigroup, polyloop) in Theorem $3.4$. This shows that second and third  nuclei of a polygroupoid (polyquasigroup, polyloop) generalize the first  nucleus of polygroupoid (polyquasigroup, polyloop). Furthermore, the fourth  nucleus of polyloop generalizes first, second and third  nuclei of polygroupoid (polyquasigroup, polyloop).

\paragraph{}
For the polygroupoid $(G, \odot)$ in Example $3.1$, the following are true:
\[N^i(G, \odot) = N^i_{\lambda}(G, \odot)= N^i_{\mu}(G, \odot) = N^i_{\rho}(G, \odot)=\emptyset~\textrm{for}~i=1,2,3,4.\]
For the polyquasigroups $(G, \odot, \nwarrow, \nearrow)$, $(G, \nwarrow, \odot, \nwarrow)$ and  $(G, \nearrow, \nearrow ,\odot)$ in Example $3.2$, the following are true:
\[N^i(G, \odot) = N^i_{\lambda}(G, \odot)= N^i_{\mu}(G, \odot) = N^i_{\rho}(G, \odot)=\emptyset~\textrm{for}~i=1,2,3,4.\]
\[N^i(G, \nwarrow) = N^i_{\lambda}(G, \nwarrow)= N^i_{\mu}(G, \nwarrow) = N^i_{\rho}(G, \nwarrow)=\emptyset~\textrm{for}~i=1,2,3,4.\]
\[N^i(G, \nearrow) = N^i_{\lambda}(G, \nearrow)= N^i_{\mu}(G, \nearrow) = N^i_{\rho}(G, \nearrow)=\emptyset~\textrm{for}~i=1,2,3,4.\]
For the polyloops $(G, \odot, \nwarrow, \nearrow,1)$, $(G, \nwarrow, \odot, \nwarrow,1)$ and  $(G, \nearrow, \nearrow ,\odot,1)$ in Example $3.3$, the following are true:
\[N^i(G, \odot) = N^i_{\lambda}(G, \odot)= N^i_{\mu}(G, \odot) = N^i_{\rho}(G, \odot)=\emptyset~\textrm{for}~i=1,2,3,4.\]
\[N^i(G, \nwarrow) = N^i_{\lambda}(G, \nwarrow)= N^i_{\mu}(G, \nwarrow) = N^i_{\rho}(G, \nwarrow)=\emptyset~\textrm{for}~i=1,2,3,4.\]
\[N^i(G, \nearrow) = N^i_{\lambda}(G, \nearrow)= N^i_{\mu}(G, \nearrow) = N^i_{\rho}(G, \nearrow)=\{1\}~\textrm{for}~i=1,2,3,4.\]
Thus, Theorem $3.1$ to Theorem $3.4$ are true for the polygroupoid, polyquasigroups and polyloops in Example $3.1$, Example $3.2$ and Example $3.3$.

\par 
\noindent
{\bf Example 3.4:} Consider $P=\{A, B, C, D, E, F, G, H, I, J, K, L\}$ with the hyperoperation '$\star$' as defined in  multiplication Table \ref{leppqua}. $(P,\star)$ is a polyquasigroup which has the Tallini 1 axiom but does not have the Tallini 2 axiom. For $(P,\star)$, the following are true: 
%\begin{center}
$N^i_{\lambda}(P,\star)=\{A, H\}$, $N^i_{\mu}(P,\star)=\{A, B,G,H\}$, $N^i_{\rho}(P,\star)=\{A, H, I, L\}$ and $N^i(P,\star)=\{A, H\}~\textrm{for}~i=1,2,3,4$.
\par 
\noindent
{\bf Example 3.5:} 	Consider $P=\{e,A, B, C, D, E, F, G, H, I, J, K, L\}$ with the hyperoperation '$\ast$' as defined in  multiplication Table \ref{lepp}. $(P,\ast)$ is a polyloop  which has the Tallini 1 axiom but does not have the Tallini 2 axiom. For $(P,\ast)$, the following are true: 
%\begin{center}
$N^i_{\lambda}(P,\ast)=\{e,A, H\}$, $N^i_{\mu}(P,\ast)=\{e,A, B,G,H\}$, $N^i_{\rho}(P,\ast)=\{e,A, H, I, L\}$ and $N^i(P,\ast)=\{e,A, H\}~\textrm{for}~i=1,2,3,4$.
\newpage
%{|l|l|l|l|l|l|l|l|l|l|l|l|l|}
%{|c|c|c|c|c|c|c|c|c|c|c|c|c|}

\begin{table}[!htbp]\caption{Multiplication table for polyquasigroup $(P,\star)$}\label{leppqua}
	\footnotesize
	\begin{tabular}{|l|l|l|l|l|l|l|l|l|l|l|l|l|}
		\hline
		$\star$  & $A$ & $B$ & $C$ & $D$ & $E$ & $F$ & $G$ & $H$ & $I$ & $J$ & $K$ & $L$ \\ 
		\hline
		$A$ & $A$ & $A,B$ & $P$ & \makecell{$A,E,$ \\ $I,D,$\\ $H,L$} & $A,E$ & \makecell{$A,E,$\\$B,F$} & \makecell{$A,G,$\\$H,B$} & $A,H$ & $A,I$ & \makecell{$A,I,$\\$B,J$} & \makecell{$A,K,$\\$L,B$} & $A,L$ \\ 
		\hline
		$B$ & $A,B$ & $B$ & \makecell{$C,G,$\\$K,B,$\\$F,J$} & $P$ & \makecell{$A,E,$\\$B,F$} & $B,F$ & $G,B$ & \makecell{$A,G,$\\$H,B$} & \makecell{$A,I,$\\$B,J$} & $B,J$ & $K,B$ & \makecell{$A,K,$\\$L,B$} \\ 
		\hline
		$C$ & $P$ & \makecell{$C,G,$\\$K,B,$\\$F,J$} & $C$ & $C,D$ & \makecell{$C,E,$\\$D,F$} & $C,F$ & $C,G$ & \makecell{$C,G,$\\$D,H$} & \makecell{$C,I,$\\$D,J$} & $C,J$ & $C,K$ & \makecell{$C,K,$\\$D,L$} \\ 
		\hline
		$D$ & \makecell{ $A,E,$ \\$I,D,$\\$H,L$} & $P$ & $C,H$ & $D$ & $E,D$ & \makecell{$C,E,$\\$D,F$} & \makecell{C$,G,$\\$D,H$} & $D,H$ & $I,D$ & \makecell{$C,I,$\\$D,J$} & \makecell{$C,K,$\\$D,L$} & $D,L$ \\ 
		\hline
		$E$ & $A,E$ & \makecell{$A,E,$\\$B,F$} & \makecell{$C,E,$\\$D,F$} & $E,D$ & $E$ & $E,F$ & $P$ & \makecell{$A,E,$\\$I,D,$\\$H,L$} & $E,I$ & \makecell{$E,I,$\\$F,J$} & \makecell{$E,K,$\\$L,F$} & $E,L$ \\ 
		\hline
		$F$ &  \makecell{$A,E,$\\$B,F$} & $B,F$ & $C,F$ & \makecell{$C,E,$\\$D,F$} & $E,F$ & $F$ & \makecell{$C,G,$\\$K,B,$\\$F,J$} & $P$ & \makecell{$E,I,$\\$F,J$} & $F,J$ & $K,F$ & \makecell{$E,K,$\\$L,F$} \\ 
		\hline
		$G$ &  \makecell{$A,G,$\\$H,B$} & $G,B$ & $C,G$ & \makecell{$C,G,$\\$D,H$} & $P$ & \makecell{$C,G,$\\$K,B,$\\$F,J$} & $G$ & $G,H$ & \makecell{$G,I,$\\$H,J$} & $G,J$ & $G,K$ & \makecell{$G,K,$\\$H,L$} \\ 
		\hline
		$H$ & $A,H$ & \makecell{$A,G,$\\$H,B$} & \makecell{$C,G,$\\$D,H$} & $H,D$ & \makecell{$A,E,$\\$I,D,$\\$H,L$} & $P$ & $H,G$ & $H$ & $I,H$ & \makecell{$G,I,$\\$H,J$} & \makecell{$G,K,$\\$H,L$} & $H,L$ \\ 
		\hline
		$I$ & $A,I$ & \makecell{$A,I,$\\$B,J$} & \makecell{$C,I,$\\$D,J$} & $I,D$ & $E,I$ & \makecell{$E,I,$\\$F,J$} & \makecell{$G,I,$\\$H,J$} & $I,H$ & $I$ & $I,J$ & $P$ & \makecell{$A,E,$\\$I,D,$\\$H,L$} \\ 
		\hline
		$J$ &  \makecell{$A,I,$\\$B,J$} & $B,J$ & $C,J$ & \makecell{$C,I,$\\$D,J$} & \makecell{$E,I,$\\$F,J$} & $F,J$ & $G,J$ & \makecell{$G,I,$\\$H,J$} & $I,J$ & $J$ & \makecell{$C,G,$\\$K,B,$\\$F,J$} & $P$ \\ 
		\hline
		$K$ &  \makecell{$A,K,$\\$L,B$} & $K,B$ & $C,K$ & \makecell{$C,K,$\\$D,L$} & \makecell{$E,K,$\\$L,F$} & $K,F$ & $G,K$ & \makecell{$G,K,$\\$H,L$} & $P$ & \makecell{$C,G,$\\$K,B,$\\$F,J$} & $K$ & $K,L$ \\ 
		\hline
		$L$ & $A,L$ & \makecell{$A,K,$\\$L,B$} & \makecell{$C,K,$\\$D,L$} & $L,D$ & $E,L$ & \makecell{$E,K,$\\$L,F$} & \makecell{$G,K,$\\$H,L$} & $H,L$ & \makecell{$A,E,$\\$I,D,$\\$H,L$} & $P$ & $K,L$ & $L$ \\ 
		\hline
	\end{tabular}
\end{table}

\begin{table}[!htbp]\caption{Multiplication table for polyloop $(P,\ast )$}\label{lepp}
	\scriptsize
	\begin{tabular}{|l|l|l|l|l|l|l|l|l|l|l|l|l|l|}
		\hline
		$\ast$ & $e$ & $A$ & $B$ & $C$ & $D$ & $E$ & $F$ & $G$ & $H$ & $I$ & $J$ & $K$ & $L$ \\ 
		\hline
		$e$ & $e$ & $A$ & $B$ & $C$ & $D$ & $E$ & $F$ & $G$ & $H$ & $I$ & $J$ & $K$ & $L$ \\ 
		\hline
		$A$ & $A$ & $A$ & $A,B$ & $P$ & \makecell{$A,E,$\\$I,D,$\\$H,L$} & $A,E$ & \makecell{$A,E,$\\$B,F$} & \makecell{$A,G,$\\$H,B$} & $A,H$ & $A,I$ & \makecell{$A,I,$\\$B,J$} & \makecell{$A,K,$\\$L,B$} & $A,L$ \\ 
		\hline
		$B$ & $B$ & $A,B$ & $B$ & \makecell{$C,G,$\\$K,B,$\\$F,J$} & $P$ & \makecell{$A,E,$\\$B,F$} & $B,F$ & $G,B$ & \makecell{$A,G,$\\$H,B$} & \makecell{$A,I,$\\$B,J$} & $B,J$ & $K,B$ & \makecell{$A,K,$\\$L,B$} \\ 
		\hline
		$C$ & $C$ & $P$ & \makecell{$C,G,$\\$K,B,$\\$F,J$} & $C$ & $C,D$ & \makecell{$C,E,$\\$D,F$} & $C,F$ & $C,G$ & \makecell{$C,G,$\\$D,H$} & \makecell{$C,I,$\\$D,J$} & $C,J$ & $C,K$ & \makecell{$C,K,$\\$D,L$} \\ 
		\hline
		$D$ & $D$ & \makecell{$A,E,$\\$I,D,$\\$H,L$} & $P$ & $C,H$ & $D$ & $E,D$ & \makecell{$C,E,$\\$D,F$} & \makecell{$C,G,$\\$D,H$} & $D,H$ & $I,D$ & \makecell{$C,I,$\\$D,J$} & \makecell{$C,K,$\\$D,L$} & $D,L$ \\ 
		\hline
		$E$ & $E$ & $A,E$ & \makecell{$A,E,$\\$B,F$} & \makecell{$C,E,$\\$D,F$} & $E,D$ & $E$ & $E,F$ & $P$ & \makecell{$A,E,$\\$I,D,$\\$H,L$} & $E,I$ & \makecell{$E,I,$\\$F,J$} & \makecell{$E,K,$\\$L,F$} & $E,L$ \\ 
		\hline
		$F$ & $F$ & \makecell{$A,E,$\\$B,F$} & $B,F$ & $C,F$ & \makecell{$C,E,$\\$D,F$} & $E,F$ & $F$ & \makecell{$C,G,$\\$K,B,$\\$F,J$} & $P$ & \makecell{$E,I,$\\$F,J$} & $F,J$ & $K,F$ & \makecell{$E,K,$\\$L,F$} \\ 
		\hline
		$G$ & $G$ & \makecell{$A,G,$\\$H,B$} & $G,B$ & $C,G$ & \makecell{$C,G,$\\$D,H$} & $P$ & \makecell{$C,G,$\\$K,B,$\\$F,J$} & $G$ & $G,H$ & \makecell{$G,I,$\\$H,J$} & $G,J$ & $G,K$ & \makecell{$G,K,$\\$H,L$} \\ 
		\hline
		$H$ & $H$ & $A,H$ & \makecell{$A,G,$\\$H,B$} & \makecell{$C,G,$\\$D,H$} & $H,D$ & \makecell{$A,E,$\\$I,D,$\\$H,L$} & $P$ & $H,G$ & $H$ & $I,H$ & \makecell{$G,I,$\\$H,J$} & \makecell{$G,K,$\\$H,L$}& $H,L$ \\ 
		\hline
		$I$ & $I$ & $A,I$ & \makecell{$A,I,$\\$B,J$} & \makecell{$C,I,$\\$D,J$} & $I,D$ & $E,I$ & \makecell{$E,I,$\\$F,J$} & \makecell{$G,I,$\\$H,J$} & $I,H$ & $I$ & $I,J$ & $P$ & \makecell{$A,E,$\\$I,D,$\\$H,L$} \\ 
		\hline
		$J$ & $J$ & \makecell{$A,I,$\\$B,J$} & $B,J$ & $C,J$ & \makecell{$C,I,$\\$D,J$} & \makecell{$E,I,$\\$F,J$} & $F,J$ & $G,J$ & \makecell{$G,I,$\\$H,J$} & $I,J$ & $J$ & \makecell{$C,G,$\\$K,B,$\\$F,J$} & $P$ \\ 
		\hline
		$K$ & $K$ & \makecell{$A,K,$\\$L,B$} & $K,B$ & $C,K$ & \makecell{$C,K,$\\$D,L$} & \makecell{$E,K,$\\$L,F$} & $K,F$ & $G,K$ & \makecell{$G,K,$\\$H,L$} & $P$ & \makecell{$C,G,$\\$K,B,$\\$F,J$} & $K$ & $K,L$ \\ 
		\hline
		$L$ & $L$ & $A,L$ & \makecell{$A,K,$\\$L,B$} & \makecell{$C,K,$\\$D,L$} & $L,D$ & $E,L$ & \makecell{$E,K,$\\$L,F$} & \makecell{$G,K,$\\$H,L$} & $H,L$ & \makecell{$A,E,$\\$I,D,$\\$H,L$} & $P$ & $K,L$ & $L$ \\ 
		\hline
	\end{tabular}
\end{table}
\newpage
%	\vspace{0.3 cm}
	\par
	{\small{
			\begin{center}{\section{CONCLUDING REMARKS}}\end{center}}}
	\vskip 0.4 cm
	\par
	\noindent %Insert your concluding remarks
The results of this work generalize existing reults in literature on nuclei of groupoid (quasigroup and loop).
	\vskip 0.5 cm
	{\small{
%		\centerline{ACKNOWLEDGEMENTS}}}
	\vskip 0.3 cm
	\noindent
	% insert your thanks if any
%	The author would like to thank the anonymous referee whose comments improved the original version of this manuscript.
	\vskip 0.5 cm
	
	\vskip 0.3 cm
	\noindent 
	
	\vskip 0.3 cm
	%\scriptsize{
		{\small{
				\centerline{REFERENCES}
		}}
		\scriptsize{
			
		}
		\par
		\noindent % Institutional address of the first author here:
		\tiny{$^{1,2}$DEPARTMENT OF MATHEMATICS,
			OBAFEMI AWOLOWO UNIVERSITY, ILE - IFE, NIGERIA
			\par 
			\noindent %E-mail address of the first author
			{\it E-mail address}: {\tt kennygilori@gmail.com} 
			\par
			\noindent % Institutional address of the second author here:
			%\scriptsize{
				$^{1}$DEPARTMENT OF MATHEMATICS,
				UNIVERSITY OF LAGOS, AKOKA, NIGERIA
				\par 
				\noindent %E-mail address of the second author
				{\it E-mail addresses}: {\tt tjayeola@oauife.edu.ng,~tgjaiyeola@unilag.edu.ng}
				\noindent % Institutional address of the second author here:
				%\scriptsize{
					\par 
				\noindent 
					$^{3}$DEPARTMENT OF MATHEMATICS AND COMPUTER SCIENCES, BRANDON UNIVERSITY, BRANDON MANITOBA, CANADA.
					\par 
					\noindent %E-mail address of the second author
					{\it E-mail addresses}: {\tt oyebolao@brandonu.ca}
					\par 
				\noindent 
				$^{4}$DEPARTMENT OF MATHEMATICS, UNIVERSITY OF IBADAN, IBADAN, NIGERIA.
				\par 
				\noindent %E-mail address of the second author
				{\it E-mail addresses}: {\tt teminiyidele@gmail.com,~ob.ogunfolu@ui.edu.ng} 
					\par 
				\noindent 
				$^{5}$DEPARTMENT OF MATHEMATICS, SCHOOL OF MATHEMATICAL SCIENCES, C. K. TEDAM UNIVERSITY OF TECHNOLOGY AND APPLIED SCIENCES, NAVRONGO, GHANA.
				\par 
				\noindent %E-mail address of the second author
				{\it E-mail addresses}: {\tt eaalhassan@cktutas.edu.gh}
			}
		\end{document}